\documentclass[a4paper,11pt]{amsart}
\usepackage{amsmath,amscd,amsfonts,amssymb,latexsym}

\newtheorem{theorem}{Theorem}
\newtheorem{proposition}[theorem]{Proposition}

\newtheorem{corollary}[theorem]{Corollary}
\newtheorem{remark}[theorem]{Remark}

\newtheorem{example}[theorem]{Example}

\def\C{\mathbb{C}}
\def\Q{\mathbb{Q}}
\def\T{\mathbb{T}}
\def\R{\mathbb{R}}
\def\Z{\mathbb{Z}}
\def\P{\mathbb{P}}
\def\qed{\hfill$\Box$\s}
\def\s{\vskip10pt}
\def\t{{\bf{t}}}
\def\map#1{\stackrel{#1}\longrightarrow}

\title[Hirzebruch class and Bia\l ynicki-Birula decomposition]
{Hirzebruch class and Bia\l ynicki-Birula decomposition}
\author{Andrzej Weber}\thanks{Supported by NCN grant 2013/08/A/ST1/00804}
\address{Department of Mathematics of Warsaw University\\
 Banacha 2, 02-097 Warszawa, Poland} \email{aweber@mimuw.edu.pl}

\begin{document}

\begin{abstract}
We establish a relation between Bia\l ynicki-Birula decomposition for $\C^*$-action and Atiyah-Bott-Berline-Vergne localization formula.
\end{abstract}

\maketitle
Suppose an algebraic torus $\C^*$ acts on a complex algebraic variety $X$.  We assume that the action is algebraic. Then a great part of the information about global invariants of $X$ is encoded in some data localized around the fixed points.
The goal of this note is to present a connection between two approaches to localization for $\C^*$-action. The homological results are related to $S^1$-action, while from $\R^*_{>0}$-action we obtain a geometric decomposition.
Algebraic (or more generally holomorphic)  actions of $\C^*\simeq S^1\times\R^*_{>0}$ have the property that the actions of $S^1$ and $\R^*_{>0}$
determine each other. From our point of view
what really matter are the invariants of fixed point set components and information about the characters of the torus acting on the normal bundles. All the data can be deduced from Bia\l ynicki-Birula decomposition of the
 fixed point sets for the finite subgroups of $\C^*$. Our main goal is to express a relation between invariants of Bia\l ynicki-Birula cell and the localized Hirzebruch class. The procedure works for smooth algebraic varieties. A part of the construction can be carried out for singular varieties.
We discuss two decompositions of Hirzebruch $\chi_y$-genus: the first one related to $S^1$-action, the second one related to $\R^*_{>0}$ flow. We show that via a limit process the second decomposition is obtained from the first one.
 That is our main result, Theorem \ref{theo1}.
 The idea to treat the generator $\t\in H^2(BS^1)$ as a real number and pass to infinity is already present in \cite{Wi}.
 A version of the main result is valid for singular varieties and in a relative context by Theorem \ref{theo2}.

\s
I would like to thank Andrzej Bia\l ynicki Birula for his valuable comments and for explaining to me the origin of the plus-decomposition. I'm also grateful to J\"org Sch\"urmann and Matthias Franz for their remarks.

\tableofcontents

\section{Homological localization associated to $S^1$--action}
The homological approach to localization was initiated by Borel \cite{Bo} in 50-ties, but a lot of ideas originate from Smith theory \cite{Sm}. Further development is due to Segal \cite{Se}, Quillen \cite{Qu}, Chang-Skjelbred \cite{ChSk}. The theory was summarized by Goresky-Kottwitz-MacPherson \cite{GKM}. This is essentially a topological theory, and the group acting is assumed to be a compact torus $T\simeq(S^1)^r$. We give the formulation by Quillen. We will discuss only cohomology with rational coefficients or later where differential forms are involved, the cohomology with real coefficients.
\begin{theorem}[\cite{Qu}, Theorem 4.4]\label{locH} Let a compact torus act on a compact topological space. Then the restriction map in equivariant cohomology \begin{equation}H^*_T(X)\to H^*_T(X^T)\end{equation}
is an isomorphism after localization in the multiplicative system generated by nontrivial characters.\label{quloc} \end{theorem}
A variant of this approach was independently developed by Atiyah-Bott \cite{AtBo1, AtBo2}, Baum-Bott \cite{BaBo}, Berline-Vergne \cite{BeVe} for manifolds. In particular there was found an inverse map of the restriction to the fixed point set. The integral of a cohomology class can be expressed as a sum of contributions concentrated on the components of the fixed point set
\begin{theorem}[Atiyah-Bott, Berline-Vergne]Let $\alpha\in H^*_T(X)$. For a component $F\subset X^T$ let $eu^T(\nu_F)\in H^*_T(F)$ denote the equivariant Euler class of the normal bundle to $F$. Then \begin{equation}\int_X\alpha=\sum_{F\text{ \rm component of }X^T}\;\int_F\frac{\alpha_{|F}}{eu^T(\nu_F)}\,.\end{equation}
\label{locth}\end{theorem}
Atiyah-Bott and also Vergne \cite{Ve} point out that the localization formula was discussed already in Witten \cite{Wi}.

All the listed papers relate the cohomological properties of the fixed point set to global invariants.
The arguments are of homological nature.
If $X$ is an algebraic variety and the action of $T$ is the restriction of an algebraic action of a complex torus $\T$, then $X^T=X^\T$. Since $H^*_T(X)\simeq H^*_\T(X)$, from the homological point of view there is no difference between $S^1$ action and $\C^*$-action.

\section{Geometric localization associated to the action of $\R^*_{>0}$}

The second, very different approach of Bia\l ynicki-Birula  is of geometric nature. The decomposition is obtained just from the action of $\R^*_{>0}\subset \C^*$, the action transverse to the $S^1$-orbits.
Let $X$ be a complete algebraic variety. Suppose that $\T=\C^*$ is the one dimensional algebraic torus. Let $\lambda\cdot x$ for $\lambda\in\C^*$ and $x\in X$ denote the action.
Let $F$ be a connected component of the fixed point set $X^\T$. Let us define
\begin{equation}X_F^+=\left\{x\in X:\lim_{\lambda\to 0} \lambda\cdot x\in F\right\}\,.\end{equation}
Note that for complete algebraic varieties the limit $\lim_{z\to 0} z\cdot x$ always exists.
To define the cell $X_F^+$ it is enough to consider the real flow generated by the fundamental field of the action. We have
\begin{equation}X_F^+=\left\{x\in X:\lim_{t\to -\infty} \exp(t)\cdot x\in F\right\}\,.\end{equation}
\begin{theorem}[Bia\l ynicki-Birula plus-decomposition, \cite{B-B1}  Theorem 4.3]\label{BBdec0}
Suppose $X$ is a complete nonsingular algebraic variety. Each set of the decomposition
\begin{equation}X=\bigsqcup_{F\text{ \rm component of }X^\T}\;X_F^+\end{equation}
is a locally closed algebraic subvariety. Taking the limit $x\mapsto\lim_{\lambda\to 0} \lambda\cdot x$ is a morphism of algebraic varieties $X_F^+\to F$.
 The map $X_F^+\to F$ is a locally trivial bundle in the Zariski topology. The fibers are isomorphic to affine spaces.
\end{theorem}
In the original formulation of Theorem \ref{BBdec0} there is an extra assumption: $X$ can be covered by $\T$-stable quasi-affine open subsets. This assumption is unnecessary by \cite[Lemma 8]{Su}.
Bia\l ynicki-Birula decomposition is valid for algebraic varieties over  an arbitrary field. Further properties  of the decomposition for projective varieties are studied in \cite{B-B3}. We will not use these properties in our considerations.

We will consider  the Bia\l ynicki-Birula decomposition for singular spaces. We assume that the singular space is contained in a smooth $\T$-variety and the action is induced from the ambient space. By \cite[Theorem 1]{Su} any quasi-projective normal $\T$-variety can be embedded equivariantly in a projective space with a linear action. If the singular space is not normal, then an embedding to a smooth $\T$-variety might not exist. For example, let $X=\P^1/\{0,\infty\}$ be the projective line with the standard $\T$-action and with the points $0$ and $\infty$ identified. Suppose $X$ can be embedded in a smooth $\T$-variety. Then by \cite[Corollary 2]{Su} there exists an affine invariant neighbourhood of each orbit. But any $\T$-invariant neighbourhood of the singular point contains the whole $X$, which is complete. Therefore $\P^1/\{0,\infty\}$ cannot be equivariantly embedded into a smooth variety.

Despite counterexamples we will consider singular spaces which might be not normal nor  projective. Instead we assume that they are equivariantly embedded in a smooth $\T$-space. This way we obtain a decomposition into locally closed subsets.

\begin{corollary}\label{BBdec} Assume that $B$ is a complete smooth $\T$-variety and $X$ is a $\T$-invariant subvariety. Then $X$ admits a decomposition into locally closed subsets
\begin{equation}X=\bigsqcup_{F\text{ \rm component of }B^\T}\;X_F^+\,,\end{equation}
where $X_F^+=X\cap B_F^+$. For each component $F\subset B$ the map $X_F^+\to F\cap X$, $x\mapsto \lim_{t\to 0}t\cdot x$ is a morphism in the category of algebraic varieties.
\end{corollary}
We can as well obtain a finer decomposition indexed by the components of $X^\T$, but we will not need that.

\section{Homological corollaries from plus-decomposition}

Suppose $X$ is smooth and complete. From Theorem \ref{BBdec} it follows that  the class of $X$ decomposes in the Grothendieck ring of varieties
\begin{equation}[X]=\sum_{F\text{ \rm component of }X^\T} [F]\cdot [\C^{n^+(F)}]\in K(Var/\C)\,,\end{equation}
where $n^+(F)$ is the dimension of the affine bundle $X_F^+\to F$.
It is equivalent to say that for any ring-valued function
\begin{equation}A:Varieties\to R\,,\end{equation} which is an additive and multiplicative invariant of algebraic varieties, we have
\begin{equation}A(X)=\sum_{F\text{ \rm component of }X^\T} A(F)\cdot A(\C)^{n^+(F)}\,.\end{equation}
In particular for Hirzebruch $\chi_y$-genus (\cite[Chapter 4]{Hi}, see the next section) we obtain
\begin{equation}\chi_y(X)=\sum_{F\text{ \rm component of }X^\T} \chi_y(F)\cdot (-y)^{n^+(F)}\,.\label{hirdec}\end{equation}
A decomposition of rational cohomology follows
when one applies virtual Poincar\'e polynomials (\cite[4.5]{Fu1}). Originally that decomposition was obtained by Bia\l ynicki-Birula by counting points of the variety reduced to a finite base field.
\begin{theorem}[\cite{B-B2}]\label{bbhom} If $X$ is smooth and complete algebraic variety, then
\begin{equation}H^*(X)=\bigoplus_{F\text{ \rm component of }X^\T} H^{*+2\,n^+(F)}(F)\end{equation}
 as graded vector spaces,
where $n^+(F)$ is the dimension of the affine bundle $X_F^+\to F$.\end{theorem}
Here again the original assumptions of \cite{B-B2} are stronger ($X$ is assumed to be projective), but the theorem holds for complete $X$ due to the Sumihiro result \cite[Lemma 8]{Su}, as already explained.
The decomposition of $X$ as an object of the category of Chow motives was proved by several authors, see the discussion in \cite[\S3]{Br}. The isomorphism of   Theorem \ref{bbhom} can be realized by a certain correspondence \cite[Theorem 3.1]{Br}
\begin{equation}X\to \bigsqcup_{F\text{ \rm component of }X^\T} F(n^+(F))\,,\end{equation} see \cite[\S16]{Fu2} for the notation.

For K\"ahler manifolds (under the assumption that $X$ is connected, $X^T\neq\emptyset$ or $H^1(X)=0$) the cohomological counterpart part of the decomposition was already known by Frankel \cite{Fr}. The decomposition of integer cohomology was proven by Carrell-Sommese \cite{CaSo}. A generalization for a certain class of singular varieties was proved by Carrell-Goresky \cite{CaGo}. We note, that for arbitrary $\C^*$-actions the decomposition is false. There are holomorphic actions on K\"ahler manifolds without fixed points (since $\C^*$ acts transitively on an elliptic curve).

\section{Relating homological and geometric decompositions}
\label{relating}
The following question arises: are the homological and geometric decompositions related? We give an answer to this question in terms of characteristic classes. Consider an invariant $A(X)$ of a smooth complete algebraic variety $X$ which is both
\begin{enumerate}\item
given by the integral of a multiplicative characteristic class of the tangent bundle,
\item additive with respect to cut and paste procedure, i.e. it factors through the Grothen\-dieck group of varieties $K(Var/\C)$.
\end{enumerate}

By \cite{BSY} the Hirzebruch $\chi_y$-genus satisfies (1) and (2).

\begin{remark}\rm The property (2) implies that:
\begin{enumerate}
\setcounter{enumi}{2}
\item for a fibration $F\hookrightarrow E\to\!\!\!\!\to B$ which is locally trivial in Zariski topology it holds
 \begin{equation}A(E)=A(F)\cdot A(B)\,.\end{equation}
\end{enumerate}
 This kind of multiplicative property was discussed in \cite{To} for complex manifolds.
By \cite[Theorem 1.2]{To} the $\chi_y$-genus is an universal invariant which is multiplicative with respect to bundles with fibers being projective spaces.
 One can check (by easy computation) that if an invariant given by a characteristic class is multiplicative with respect to algebraic $\P^1$- and $\P^2$-bundles over $\P^n$'s, then such an invariant has to be
(up to a change of variables) the Hirzebruch $\chi_y$-genus. It follows that the $\chi_y$-genus is an universal invariant satisfying (1) and (3) for algebraic varieties. Since (3) is weaker than (2), the $\chi_y$-genus is also an universal  invariant satisfying (1) and (2).
\end{remark}

The $\chi_y$-genus and the Hirzebruch class $td_y$ provide a bridge between geometric decompositions of varieties and decompositions on the cohomology level. For $y=0$ we obtain the  Todd class which was studied by Brion and Vergne in \cite{BrVe} from the equivariant point of view in the case of toric varieties.
\s
To describe the simplest case of our main result assume that $X$ is smooth, complete and the fixed point set is finite. Then for $p\in X^\T$ the associated subset $X_{\{p\}}^+$  is an algebraic cell isomorphic to an affine space. The cell dimension $n^+(p)$ can be expressed in terms of the tangent representation $T_pX$. Let
 \begin{equation}T_pX=\bigoplus_{w\in \Z} V_w\end{equation}
 be the decomposition into weight subspaces of the action. Then \begin{equation}n^+(p)=\sum_{w>0}\dim(V_w)\,.\end{equation}
Let \begin{equation}P_X(t)=\sum_{i=0}^{2\,\dim(X)}\dim(H^i(X;\Q))\,t^i\end{equation}
be the Poincar\'e polynomial of $X$. When $X$ admits a decomposition into algebraic cells
then $P_X$ (after the substitution ${y=-t^2}$) is equal to the Hirzebruch $\chi_y$-genus \cite[\S IV.15.5]{Hi}
\begin{equation}P_X(t)=\chi_y(X)=\sum_{p=0}^{\dim(X)}\chi(X;\Omega^p_X)y^p~_{|y=-t^2}\,.\end{equation}
Each fixed point of the $\T$-action gives rise to a generator of cohomology and contributes to $\chi_y$-genus:
\begin{equation}\chi_y(X)=\sum_{p\in X^\T}(-y)^{n^+(p)}\,.\end{equation}

On the other hand the $\chi_y$-genus can be computed by integration of the Hirzebruch class $td_y(X)$.
We recall, that the Hirzebruch class of a complex manifold is defined as a characteristic class of the tangent bundle
\begin{equation}td_y(X):=td_y(TX)\,.\end{equation}
The (nonreduced) class $td_y(-)$ is the multiplicative characteristic class associated to the power series
\begin{equation}x\frac{1+y\,e^{-x}}{1-e^{-x}}=x\left(\frac{1+y}{1-e^{-x}}- y\right)\,,\end{equation}
see \cite[\S5.1]{Hy} for a short presentation.
If the torus acts on $X$ then naturally the Hirzebruch class (as any characteristic class) lifts to equivariant cohomology. The properties of the equivariant Hirzebruch class are studied in \cite{We}. An example of computation is given in \cite{MiWe}.
Having the equivariant version of the Hirzebruch class $td^\T_y(X)$ we can apply Atiyah-Bott-Berline-Vergne formula. Strictly speaking we have a formula for push forward $p_*:H^*_\T(X)[y]\to H^{*-2\dim X}_\T(pt)[y]$.
By \cite{Mu} the class $p_*td^\T_y(X)$ is concentrated in degree zero and it is equal to $\chi_y(X)$. Thus we have
\begin{equation}\chi_y(X)=\int_X td^\T_y(X)=\sum_{p\in X^\T}\frac{td^\T_y(X)_{|p}}{eu^\T(T_pX)}\,.\end{equation}
The contribution to the $\chi_y$ genus corresponding to a point $p\in X^\T$ is equal to
\begin{equation}\frac{td^\T_y(X)_{|p}}{eu^\T(T_pX)}=\prod_{w\in \Z}\left(\frac{1+y\,e^{-w\t}}{1-e^{-w\t}}\right)^{\dim(V_w)}\,,\end{equation}
where $\t$ is the generator of $H^*_\T(pt)\simeq\R[\t]$ (we will consider cohomology with real coefficients). The distinguished generator $\t\in H^2(B\T)$ corresponds to the identity character $\T=\C^*\to\C^*$.
Observe
\begin{equation}\lim_{t\to -\infty}\frac{1+y\, e^{-wt}}{1- e^{-wt}}=\left\{\begin{matrix}-y&\text{for } w>0\\1&\text{for } w<0\end{matrix}\right.\,.\end{equation}
Let us treat $\t$ as a real number. Then
\begin{equation}\lim_{\t\to-\infty} \frac{td^\T_y(X)_{|p}}{eu^\T(T_pX)}=(-y)^{n^+(p)}\,,\end{equation}
which is exactly the contribution to $\chi_y(X)$ coming from the Bia\l ynicki-Birula decomposition. Therefore we can loosely say:
\begin{corollary} The Bia\l ynicki-Birula decomposition of $\chi_y$-genus is the limit of the Atiyah-Bott-Berline-Vergne localization formula.\end{corollary}

\begin{example}\rm \label{exp1} Let $X=\P^1=\C\cup\{\infty\}$ with the standard action of $\C^*$.
Then $X^\T=\{0,\infty\}$ and
\begin{equation} X_{\{0\}}=\C\,,\quad X_{\{\infty\}}=\{\infty\}\,. \end{equation}
We obtain a decomposition in $K(Var/\C)$
\begin{equation}\tag{{\bf motivic}}\qquad[\P^1]=[\C]+[pt]\,. \end{equation}
The resulting decomposition of $\chi_y$ is the following:
\begin{equation}\tag{{\bf BB}}\qquad\chi_y(\P^1)=(-y)+1\,. \end{equation}
The Atiyah-Bott or Berline-Vergne integration formula allows to decompose
\begin{equation}\tag{{\bf ABBV}}\qquad\chi_y(\P^1)=\frac{1+y e^{-\t}}{1- e^{-\t}}+\frac{1+y e^{\t}}{1- e^{\t}}\,.\end{equation}
Allowing $\t\to-\infty$ we see that {\bf (ABBV)} decomposition converges to {\bf (BB)} decomposition.
\end{example}

\section{Relative Hirzebruch class}\label{Kvar}

 Characteristic classes of equivariant vector bundles appear in literature quite often, see e.g. \cite[\S7]{BGV} where a differential definition is given. For a compact connected Lie group $K$ acting via symplectic transformations on a symplectic manifold $X$
the equivariant Todd class $td^K(X)=td^K(TX)\in \prod_{i=0}^\infty H^i_K(X)$  was studied in \cite[\S2]{JeKi}.
In the formulation of our main theorem, Theorem \ref{theo1}, the Hirzebruch class $td^\T_y$ for singular varieties appears.  The non-equivariant case was studied in \cite{BSY}.
The equivariant version for a torus action was developed in \cite{We}. Let us list the formal properties which determine
 this class. Let us assume that $B$ is a  complex smooth algebraic $\T$-variety. For any equivariant map of $\T$-varieties $f:X\to B$ we have a cohomology class \begin{equation}td_y^\T(f:X\to B)\in \hat H^*_T(B)[y]:=\prod_{k=0}^{+\infty}H^k_\T(B)[y]\end{equation}
 satisfying:
 \begin{itemize}
 \item[(i)] if $X=B$ and $f=id_X$, then $td_y^\T(f:X\to B)=td_y^\T(TX)$ is defined by the equivariant characteristic  class associated to the power series $x\frac{1+y\,e^{-x}}{1-e^{-x}}$.

\item[(ii)] if $g:B_1\to B_2$ and $f:X\to B_1$ are equivariant maps and $g$ is proper, then \begin{equation}td_y^\T(g\circ f:X\to B_2)=g_*td_y^\T(f:X\to B_1)\,,\end{equation}

\item[(iii)] if $U$ is a $\T$-stable open set of $X$, $Y=X\setminus U$, then
\begin{equation}td_y^\T(f:X\to B)=td_y^\T(f_{|U}:U\to B)+td_y^\T(f_{|Y}:Y\to B)\,.\end{equation}
\end{itemize}
The above properties  indicate how to calculate $td_y^\T(f:X\to B)$.
The Hirzebruch class of $f:X\to B$ is computed inductively with respect to the dimension of $X$.
We may decompose $X$ into smooth $\T$-invariant strata and use the additivity property (iii). Therefore it is enough to assume that $X$ is smooth. We claim that there exists a smooth $\T$-variety $\overline X$ containing $X$ and an equivariant extension $\overline{f}:\overline{X}\to B$ of $f$ which is  proper.
 To construct $\overline X$ we find (by  \cite[Theorem 3]{Su}) an equivariant completion $Y_1\supset X$. Let $Y_2$ be the closure of the graph $f:X\to B$ in $Y_1\times B$. The variety $Y_2$ contains an open subset isomorphic to $X$ and the restriction of the projection to $B$ is a proper extension of $f$. Now we resolve the singularities of $Y_2$ in an equivariant way due to \cite{BiMi} and obtain the desired variety $\overline X$.
Now \begin{equation}td_y^\T(f:X\to B)=td_y^\T(\overline{f}:\overline {X}\to B)-td_y^\T(\overline{f}_{|\partial X}:\partial X\to B)\,,\end{equation}
where $\partial X=\overline {X}\setminus X$. The first summand can be computed from the property (i) and (ii). Since $\dim(\partial X)<\dim(X)$ the term $td_y^\T(\overline{f}_{|\partial X}:\partial X\to B)$ is computed by the inductive assumption. Of course it is a nontrivial fact that the obtained result does not depend on the choices made.
\begin{remark}\rm By the equivariant Chow lemma \cite[Theorem 2]{Su} we can modify $\overline X$ and assume that it is quasi-projective.\end{remark}

The procedure described above  can be purely formally expressed as follows.
Consider the Grothendieck group $K^\T(Var/B)$ of varieties over $B$ with a compatible torus action. This group is ge\-ne\-rated by the classes of equivariant maps $[f:X\to B]$. The generators satisfy the relation: for any open $\T$-stable subset $U\subset X$, $Y=X\setminus U$: \begin{equation}[f:X\to B]=[f_{|U}:U\to B]+[f_{|Y}:Y\to B]\,.\end{equation}
By  equivariant completion and equivariant resolution  the group $K^\T(Var/B)$ is gene\-rated by the proper maps $[X\to B]$, where $X$ is smooth. We do not go into the structure of the group (in fact a ring)
$K^\T(Var/B)$. We treat it here just as an auxiliary formal construction\footnote{The generators and relations of non-equivariant $K(Var/B)$ were described in \cite{Bi}.
 We know that $K^\T(Var/B)$ is generated by the proper equivariant maps from smooth, quasi-projective  varieties.
 A description of relations is not present in the literature.  The relations among generators are irrelevant for us.}.
\s
The main property of the Hirzebruch class is that the assignment
$$\begin{matrix}\big\{ \text{proper equivariant maps from smooth varieties to } B\big\}&~~\longrightarrow~~&\hat H^*_\T(B)[y]\\
f:X\to B & ~~\longmapsto~~ & f_*(td^\T_y(X))\end{matrix}$$
extends to a homomorphism of groups
\begin{equation}K^\T(Var/B)\to \hat H^*_\T(B)[y]\,.\end{equation}
In the non-equivariant case this is the main theorem of \cite{BSY}. The equivariant case follows automatically. It is discussed in \cite{We}.

\section{The main result}
We will generalize the simple calculation of the Example \ref{exp1} to the case of possibly singular varieties and actions with arbitrary fixed point set.
\begin{theorem} \label{theo1}Let $\T=\C^*$ and let $B$ be a smooth complete algebraic variety with $\T$-action, $X\subset B$ a possibly singular closed invariant subvariety. Let $F\subset B^\T$ be a component of the fixed point set.
 Then the limit of the localized Hirzebruch class is equal to the Hirzebruch class of $X_F^+=B_F^+\cap X$
\begin{equation}\lim_{\t\to-\infty}\frac{td^\T_y(X\to B)_{|F}}{eu^\T(\nu_F)}=td_y(X^+_{F}\to F)\in H^*(F)[y]
\,.\end{equation}
\end{theorem}

The substitution of the cohomology class $\t\in H^2(B\T)$ by a real number needs some justification.
This is carefully explained by Propositions \ref{zbieznosc}-\ref{niezerowosc} and Theorem \ref{rozne}.

\begin{remark}\rm Regardless of the formal proof we would like to note that the convention to treat $\t$ as a {\it real number}, not as a {\it cohomology class},
is present in many places, see for example Witten \cite[\S3]{Wi}. Also the idea to look at the limit when $\t\to -
\infty$ is present in Witten's paper.
\end{remark}

Since the Hirzebruch class $td_y^\T(X\to B)$ is constructed via resolution of $X$ we are led to prove a slightly more general functorial formulation. In Theorem \ref{theo2} we consider not a subvariety, but an arbitrary map $X\to B$. This extension is forced by  the proof of the theorem based on the properties (i)-(iii) of \S \ref{Kvar}.
By these properties it is enough to consider the case when $X=B$, and the result follows from the calculus of characteristic classes.

We summarize our result in following way:
\begin{corollary}\label{theofunc} Let $B$ be a smooth and complete algebraic variety with a $\T$-action. Denote by $K^\T(Var/B)$ the Grothendieck group of varieties over $B$ equipped with a compatible torus action. Let $F\subset B^{\T}$ be a component of the fixed point set.
The following diagram relating Grothendieck groups with cohomology is commutative:
\begin{equation}\begin{matrix} &{\beta_F}&\\
K^\T(Var/B)&\longrightarrow&K(Var/F)\\ \\
{td^\T_y}\downarrow\phantom{{td^\T_y}}&&\phantom{{td_y}}\downarrow{td_y}\\ \\
\widehat H^*_\T(B)[y]&\dashrightarrow&H^*(F)[y],\\
&{\lim_F}\end{matrix}\end{equation}
where \begin{equation}\beta_F([f:X\to B])=[f_{|f^{-1}(B^+_F)}:f^{-1}(B^+_F)\to F]\end{equation} and the map ${\lim}_F$ given by
\begin{equation}{\lim}_F(\alpha)=\lim_{\t\to-\infty}\frac{\alpha_{|F}}{eu^\T(\nu_F)}\end{equation} is well defined for $\alpha$ belonging to the image of ${td^\T_y}$.
\end{corollary}

For a singular variety $X$ the genus $\chi_y(X)$ is by definition equal to $td_y(X\to pt)\in H_*(pt)[y]=\R[y]$. (In fact $td_y(X\to pt)$ has integer coefficients.)
Having Theorem \ref{theo1} for granted we conclude that Atiyah-Bott-Berline-Vergne decomposition of $\chi_y$-genus
\begin{equation}\chi_y(X)
=\sum_{F\text{ \rm component of }B^T}\int_F\frac{td^\T_y(X\to B)_{|F}}{eu^\T(\nu_F)}\end{equation}
converges to the Bia\l ynicki-Birula decomposition
\begin{equation}\chi_y(X)=\sum_{F\text{ \rm component of }B^T}\chi_y(B^+_{F}\cap X)\end{equation}
as $\t$ tends to $-\infty$.

\begin{remark}\rm
Application of the Laurent series $\frac{1+y\,e^{-\t}}{1-e^{-\t}}$ might be misleading for the reader. In fact by no means
we rely on analytic properties of that series. It just serves as a bridge from equivariant $K$-theory to equivariant cohomology. The series $e^\t\in \hat H^*_\T(pt)$ is equal to the Chern character of the natural representation of $\T$. In other words it is the image of a distinguished generator of equivariant $K$-theory under the equivariant Chern character map\footnote{The equivariant $K$-theory $K_\T(pt)\simeq Rep(\T)$ of \cite{Se} should not be confused with $K$-theory of varieties $K^\T(Var/pt)=K^\T(Var/\C)$.} \begin{equation}ch:K_\T(pt)\to \hat H^*_\T(pt)\,.\end{equation}  A better environment
for our theorem would be equivariant $K$-theory of $B$. Then, instead of $e^\t$, we might use a variable $\theta$ and all the expressions
involved would be rational functions in $\theta$.
Nevertheless we have decided to work with cohomology since we believe that this setup is wider spread in literature.\end{remark}

In the section \S\ref{podgrupy} we study the Hirzebruch class of the fixed points $X^G$ for a subgroup $G\subset \T$. We show that for a smooth variety $X$ it is possible to read the $\chi_y$-genera of the cells
$(X^G)^+_F$
from the localized equivariant class.

\section{Asymptotic of the characteristic class $ch^\T\Lambda_y$}

Assume that $X$ is a smooth $\T$-variety. Our goal is to prove that \begin{equation}\frac{td^\T_y(X\to B)_{|F}}{eu^\T(\nu_F)}\end{equation} makes sense after the substitution $\t=t_0\in \R\setminus\{0\}$. Then we compute the  limit when $t_0\to-\infty$. The content of this paragraph can be summarized as follows: by the splitting principle we assume that the normal bundle $\nu_F$ decomposes into equivariant line bundles $\bigoplus_i\xi_i$. Let $x_i=c_1(L_i)\in H^2(F)$ and let $w_i\in \Z$ be the weight of $\T$-action on $\xi_i$.
 Then
\begin{equation}\frac{td^\T_y(X)_{|F}}{eu^\T(\nu_F)}=td_y(F)\prod_i\frac{1+y\,e^{-w_i\t}e^{-x_i}}{1-e^{-w_i\t}e^{-x_i}}\in H^*(F)[[\t]][\t^{-1},y]\,.\end{equation}
We see that the factors of the product after substitution $\t=t_0\in\R\setminus\{0\}$ become \begin{equation}\frac{1+y\,e^{-w_it_0}e^{-x_i}}{1-e^{-w_it_0}e^{-x_i}}\in H^*(F)[y]\end{equation} and converge when $t_0\to-\infty$ to $-y$ or 1, depending on the sign of $w_i$. To be fully precise we divide the proof into several steps.\s

 Let $\xi$ be a complex vector bundle over a paracompact topological space $F$ of finite dimension. Denote by $\Lambda_y(\xi)$ the formal combination
 \begin{equation}\sum_{p=0}^{\dim( \xi)}\Lambda^p(\xi)\,y^p\end{equation}
 considered as an element of the topological $K$-theory extended by a variable $y$, i.e.~an element of $K(F)[y]=K(F)\otimes \Z[y]$.
 The operation $\Lambda_y$ is of exponential type. It converts direct sums to tensor products:
 \begin{equation}\Lambda_y(\xi\oplus \eta)=\Lambda_y(\xi)\otimes \Lambda_y( \eta).\end{equation}
 Applying Chern character to $\Lambda_y(\xi)$ we obtain a characteristic class \begin{equation}ch\Lambda_y(\xi)\in H^*(F)[y].\end{equation}
 By the exponential property of $\Lambda_y$ we have
\begin{equation}ch\Lambda_y(\xi\oplus \eta)=ch\Lambda_y(\xi)\cdot ch\Lambda_y( \eta)\,.\end{equation}
\s
 Let $\xi$ be a $\T$-equivariant vector bundle over a base $F$. We will apply the localization Theorem \ref{locth}, and thus we concentrate on analysis of equivariant vector bundles over a base $F$  with trivial action. Every equivariant bundle over such base space decomposes into a direct sum of subbundles corresponding to various characters of $\T$
 \begin{equation}\xi=\bigoplus_{w\in \Z}\xi^w\otimes \C_w\,.\end{equation}
 Here we assume that the action of $\T$ on $\xi^w$ is trivial.
 The line bundle $\C_w$ is trivial as a vector bundle with the action of $z\in\T$  via multiplication by $z^w$, for the given $w\in\Z$.

 Recall that we have chosen a generator of the cohomology ring $H^*(B\T)=\R[\t]$  (the generator $\t\in H^2(B\T)$ corresponding to the identity character $\T=\C^*\to\C^*$).
 Our goal is to describe the asymptotic behavior of the equivariant characteristic class $ch^\T\Lambda_y$ which belongs to
 \begin{equation}ch^\T\Lambda_y(\xi)\in \widehat H^*_\T( F)[y]=\widehat H^*(B\T\times F)[y]=H^*(F)[[\t]][y]\,.\end{equation}
 We use the completed cohomology $\widehat H^*=\prod_{k=0}^{+\infty} H^k$, since the classifying space is of infinite dimension.
 We will specialize the cohomology class $\t$ treated as a formal variable to a real number $t_0$.
 First of all we have to show that such a specialization makes sense.
\begin{proposition}\label{zbieznosc} The formal power series $ch^\T\Lambda_y(\xi)H^*(F)[[\t]][y]$ converges for all substitutions $\t=t_0\in \R$.\end{proposition}

The resulting specialization will be denoted by $ch^{t_0}\Lambda_y(\xi)\in H^*(F)[y]$.
For $y=-1$ we will need to invert the classes $ch^\T\Lambda_{-1}(\xi)$ and $ch^{t_0}\Lambda_{-1}(\xi)$

\begin{proposition}\label{niezerowosc}  Assume that all weights appearing in $\xi$ are nonzero. Then the class $ch^\T\Lambda_{-1}(\xi)$ is invertible in $H^*(F)[[\t]][\t^{-1}]$. Moreover, specializing
 $\t$ to a real number $t_0\not=0$ we obtain an invertible element in $ch^{t_0}\Lambda_{-1}(\xi)\in H^*(F)$.\end{proposition}

Next, we will let $t_0$ tend to $-\infty$.

\begin{theorem}\label{rozne}Let
\begin{equation}n^+=\sum_{w>0}\dim(\xi^w)\,,\qquad n^-=\sum_{w<0}\dim(\xi^w).\end{equation}
 and assume that all weights appearing in $\xi$ are nonzero.
Then there exist the following limits in $ H^*(F)[y]$:
\begin{equation}\lim_{t_0\to \infty}\frac{ch^{t_0}\Lambda_{y}(\xi)}{ch^{t_0}\Lambda_{-1}(\xi)}=(-y)^{n^+}\,,\end{equation}
\begin{equation}\lim_{t_0\to -\infty}\frac{ch^{t_0}\Lambda_{y}(\xi)}{ch^{t_0}\Lambda_{-1}(\xi)}=(-y)^{n^-}\,.\end{equation}
In particular the limits belong to $H^0(F)[y]$.\end{theorem}

\section{Proofs of asymptotic properties of $ch^\T\Lambda_y$ }

{\it Proof of Proposition \ref{zbieznosc}.}
 First assume, that $\xi=\xi^w\otimes \C_w$.  We have
\begin{align}ch^\T\Lambda_y(\xi\otimes \C_w)&=
\sum_{p=0}^{\dim (\xi)} ch^\T(\Lambda^p \xi^w\otimes\C_{pw})y^p\\ &=
\sum_{p=0}^{\dim (\xi)} ch(\Lambda^p \xi^w)\exp(pw\t)y^p\in H^*(F)[[\t]][y]\end{align} and since the power series $\exp(\t)$ converges for any argument, $ch^\T\Lambda_y(\xi\otimes \C_w)$ specializes to
\begin{equation}\sum_{p=0}^{\dim (\xi)} ch(\Lambda^p \xi)
e^{p w t_0}y^p\in H^*(F)[y]\,.\end{equation}
The claim for an arbitrary bundle is obtained by the exponential property of $ch^\T\Lambda_y$.
\qed

{\it Proof of Proposition \ref{niezerowosc}.} It is enough to show the invertibility for $\xi=\xi^w\otimes \C_w$.  Since we assume that the base is paracompact, there exists a bundle $\eta$ over $F$, such that $\xi^w\oplus \eta$ is a trivial bundle of some rank $N\in \Z$.
Then \begin{equation}(\xi^w\otimes \C_w)\oplus (\eta\otimes \C_w)\simeq \C_w^{\oplus N},\end{equation} thus
\begin{equation}ch^\T\Lambda_{-1}(\xi^w\otimes \C_w)\cdot ch^\T\Lambda_{-1} (\eta\otimes \C_w)
=ch^\T\Lambda_{-1}( \C_w^{\oplus N})=(1-\exp(w\t))^N\,.\end{equation}
The power series $1-\exp(w\t)=-(w\t+\frac{(w\t)^2}2+\dots)$ is invertible, since $w\not=0$ and we invert $\t$.
Therefore \begin{equation}(ch^\T\Lambda_{-1}(\xi\otimes \C_w))^{-1}=(1-\exp(w\t))^{-N}\cdot ch^\T\Lambda_{-1}(\eta\otimes \C_w)\in H^*(F)[[\t]][\t^{-1}]\,.\end{equation}
Specializing $\t= t_0\not=0$ we have $1-\exp(w\t)\mapsto 1-e^{wt_0}\not=0\in\R$ and we obtain invertibility of $ch^{t_0}\Lambda_{-1}(\xi)\in H^*(F)$.
\qed

{\it Proof of Theorem \ref{rozne}.} We will concentrate on the first limit. The second one is computed by changing $z\in\T $ to $z^{-1}$ .
By the exponential property of $ch^{t_0}\Lambda$ it is enough to assume that $\xi=\xi^w\otimes \C_w$.
 Moreover, by the splitting principle we can assume, that $\dim(\xi)=1$.
Precisely, by \cite[\S17.5, Prop. 5.2]{Hu} there exists a space $F'$ and a map $f:F'\to F$, which induces a monomorphism on cohomology, and the bundle $f^*\xi$ decomposes into linear summands. By the exponential property of $ch^{t_0}\Lambda_y(-)$  it is enough to compute the limit (40)  in $H^*(F')$ for each summand of $f^*\xi$. Therefore we can assume that  $\dim\xi=1$. We note that
\begin{equation}\frac{ch^{t_0}\Lambda_{y}(\xi)}{ch^{t_0}\Lambda_{-1}(\xi)}=\frac{ch^{t_0}\Lambda_{y}(\xi)}{ch^{t_0}\Lambda_{-1}(\C_w)}
\cdot \left(\frac{ch^{t_0}\Lambda_{-1}(\xi)}{ch^{t_0}\Lambda_{-1}(\C_w)}\right)^{-1}.\label{product}\end{equation}
 We will show

\begin{equation}
 \lim_{t_0\to\infty}\frac{ch^{t_0}\Lambda_y(\xi)}{ch^{t_0}\Lambda_{-1}(\C_w)}=
 \left\{\begin{matrix}-y\,ch(\xi)\quad &\text{if } w>0\hfill\\
1\quad &\text{if } w<0.\end{matrix}\right.
\label{product1}
\end{equation}
Indeed, set $c_1(\xi)=x$. Then
\begin{equation}\frac{ch^{t_0}\Lambda_y(\xi)}{ch^{t_0}\Lambda_{-1}(\C_w)}=\frac{1+y e^{t_0w}\exp(x)}{1- e^{t_0w}}\,.\end{equation}
The denominator is a real number and we obtain the limits as claimed.
Substituting $y=-1$ we obtain the formula for the second factor of (\ref{product1}):
\begin{equation}
 \lim_{t_0\to\infty}\frac{ch^{t_0}\Lambda_{-1}(\xi)}{ch^{t_0}\Lambda_{-1}(\C_w)}=
 \left\{\begin{matrix}ch(\xi)\quad &\text{if } w>0\hfill\\
1\quad &\text{if } w<0.\end{matrix}\right.
\label{product2}\end{equation}
Taking the quotient of the limits (\ref{product1}) and (\ref{product2}) we prove the theorem.
\qed
\section{Asymptotic of the Hirzebruch class}
Let us now
assume that $w=0$ does not appear among the weights.
The equivariant Hirzebruch class of a vector bundle $\xi$ can be expressed by the class $ch^\T\Lambda_y(\xi^*)$ and the Euler class $eu^\T(\xi)$
\begin{equation}td^\T_y(\xi)=td^\T(\xi)ch^T\Lambda_y(\xi)=eu^\T(\xi)\frac{ch^\T\Lambda_y(\xi^*)}{ch^\T\Lambda_{-1}(\xi^*)}\in H^*(F)[[\t]][\t^{-1},y]\,.\end{equation}
The Euler class is invertible in the equivariant cohomology localized in $\t$ (the same proof as of Proposition \ref{niezerowosc}, see also \cite{EdGr}) and we can write
the quotient
 class
 \begin{equation}\frac{td^\T_y(\xi)}{eu^\T(\xi)}
=\frac{ch^\T\Lambda_y(\xi^*)}{ch^\T\Lambda_{-1}(\xi^*)}\in H^*(F)[[\t]][\t^{-1},y]\end{equation} expressed only by the class $ch^\T\Lambda_y(\xi^*)$.
From Theorem \ref{rozne} we obtain
\begin{corollary}\label{granice}Let $\xi$ be an equivariant  vector bundle without  eigenvectors of weight $0$.
The following limits in $H^*(F)[y]$ exist and
\begin{align}
\lim_{t_0\to \infty}\left(\frac{td^\T_y(\xi)}{eu^\T(\xi)}\right)_{|\t=t_0}&=(-y)^{n^-},\\
\lim_{t_0\to -\infty}\left(\frac{td^\T_y(\xi)}{eu^\T(\xi)}\right)_{|\t=t_0}&=(-y)^{n^+}\,.\end{align}
In particular the limits belong to $H^0(F)[y]$.
\end{corollary}
Note that
the roles of $\infty$ and $-\infty$ are interchanged since in the definition of $td_y$ the dual bundle $\xi^*$ appears.

\section{The limit of the localized Hirzebruch class - smooth case}
Let $X$ be a complex manifold with a holomorphic action of $\T$. Let $F\subset X^\T$ be a component of the fixed point set. The normal bundle $\nu_F$ decomposes into the eigen-subbundles of the torus action
\begin{equation}\nu_F=\bigoplus_{w\not=0}\nu_F^w\otimes \C_w\,.\end{equation}
The equivariant Hirzebruch class of $X$ restricted to $F$ is equal to the cup-product
\begin{equation}td^\T_y(X)_{|F}=td^\T_y(F)\cdot td_y^\T(\nu_F)\,.\end{equation}
The contribution appearing in the Berline-Vergne or Atiyah-Bott formula corresponding to the component $F$  is equal to
\begin{equation}\left(\frac{td^\T_y(X)_{|F}}{eu^\T(\nu_F)}\right)_{\t= t_0}=td_y(F)\cdot \left(\frac{td_y^\T(\nu_F)}{eu^\T(\nu_F)}\right)_{\t= t_0}
\end{equation}
By Corollary \ref{granice} we conclude:

\begin{theorem}\label{gladloc} Let $X$ and $F$ be as above, and let \begin{equation}n^+(F)=\sum_{w>0}\dim\left((\nu_F)_w\right)\,,\qquad n^-(F)=\sum_{w<0}\dim\left((\nu_F)_w\right)\end{equation} be the dimensions of the subbundles of $\nu_F$ corresponding to the positive and negative weights. Then
\begin{equation}\lim_{t_0 \to -\infty}\left(\frac{td_y^\T(X)_{|F}}{eu^\T(\nu_F)}\right)_{\t= t_0 }=
(-y)^{n^+(F)}\,td_y(F)\in H^*(F)[y]\,,\end{equation}
\begin{equation}\lim_{t_0 \to \infty}\left(\frac{td_y^\T(X)_{|F}}{eu^\T(\nu_F)}\right)_{\t= t_0 }=
(-y)^{n^-(F)}\,td_y(F)\in H^*(F)[y]\,.\end{equation}\end{theorem}

\section{ Relative localization for singular varieties}
This section is devoted to the proof of the main result, Theorem \ref{theo1} and its refined version, Theorem \ref{theo2}.

We fix $F\subset B^\T$, a component of the fixed point set. Denote by $i_F:B^+_F\to B$ the inclusion and $r_F: B_F^+\to F$ the retraction induced by the action of $\T$.
Let \begin{equation}f_F=r_F\circ f_{|f^{-1}(B_F^+)}:f^{-1}(B_F^+)\to F\end{equation} be the composition of the map $f$ with the retraction.
We consider two classes in $H_*(F)[y]$:\s
\begin{description}
\item [\bf 1) $\Phi_F$]
 The first characteristic class is the limit
 \begin{equation}\Phi_F(X\to B)=\lim_{t\to -\infty}\left(\frac{td^\T_y(X\to B)_{|F}}{eu^\T(\nu_F)}\right)_{\t= t_0}\,.\end{equation}
The limit exists for smooth and complete $X$ by Theorem \ref{rozne}.
Since the Hirzebruch class satisfies the additivity condition  (\S\ref{relating}, iii)  the limit defines a homomorphism
 $K^\T(Var/B)\to H_*(F)[y]$.
In the notation of Corollary \ref{theofunc} we have $\Phi_F=\lim_F\circ\, td_y^\T$.
\s\item[\bf 2) $\Psi_F$]
 The second class is equal to
 \begin{equation}\Psi_F(f:X\to B)=td_y(f_F:f^{-1}(B_F^+)\to F)\,.\end{equation}
This characteristic class is the composition
\begin{equation}K^\T(Var/B)\map{} K(Var/B)\map{{(i_F)}^*} K(Var/B^+_F)\map{{(r_F)}_*} K(Var/F)\map{td_y} H_*(F)[y]\,.\end{equation}
In the notation of Corollary \ref{theofunc} we have $\Psi_F= td_y\circ \beta_F $.
\s\end{description}
Both characteristic classes give rise to group homomorphisms
\begin{equation}\Phi_F,\Psi_F:K^\T(Var/B)\map{}H_*(F)[y]\,.\end{equation}
Our main result is the following:
\begin{theorem} \label{theo2} For any equivariant  map $f:X\to B$ the classes $\Phi_F(X\to B)$ and
$\Psi_F(X\to B)$
are equal:
\begin{equation}td_y(f_F:f^{-1}(B_F^+)\to F)=\lim_{t_0\to -\infty}\left(\frac{td^\T_y(f:X\to B)_{|F}}{eu^\T(\nu_F)}\right)_{\t= t_0}\,.\end{equation}
\end{theorem}

\proof 
As explained in \S\ref{Kvar} the group $K^\T(Var/B)$ is generated by the proper maps from smooth manifolds. Since we assume that $B$ is complete it is enough to check the equality for equivariant maps to $B$ from  smooth compete varieties.
Due the obvious property
\begin{equation}f^{-1}(B_F^+)=\bigsqcup_{F'\;\text{component of}\;f^{-1}(F)\cap X^\T} X_{F'}^+\,.\label{additivity}\end{equation}
and functoriality of the class $td_y^\T$ with respect to the push-forwards (property (ii) of \S\ref{Kvar}) one can assume that $B=X$, $f=Id_X$. The retraction $r_F:X_F^+\to F$ is an affine bundle with the fiber isomorphic to $\C^{n^+(F)}$, by Theorem \ref{BBdec}. Therefore \begin{equation}\Psi_F(id:X\to X)=td_y(r_F:X_F^+\to F)=(-y)^{n^+(F)}td_y(F)\,.\end{equation} This can be seen already on the level of the Grothendieck groups. We have \begin{equation}[X_F^+\to F]=[id:F\to F]\times [\C^{n^+(F)}\to pt]\,.\end{equation} (Here $\times:K(Var/F)\times K(Var/\{pt\})\to K(Var/(F\times \{pt\}))=K(Var/F)$ is given by the Cartesian product of maps.) On the other hand by Theorem \ref{gladloc}
\begin{equation}\Phi_F(id:X\to X)=\lim_{t_0 \to -\infty}\left(\frac{td_y(X)_{|F}}{eu^\T(\nu_F)}\right)_{\t= t_0 }=
(-y)^{n^+(F)}\,td_y(F)\,.\end{equation}\qed

\begin{remark}\rm
Taking the limit when $t_0\to+\infty$ one obtains the Hirzebruch class of the map
\begin{equation}f^{-1}(B^-_F)\to F\end{equation}
corresponding to the opposite minus-decomposition of $B$
\begin{equation}B_F^-=\left\{b\in B: \lim_{\lambda\to \infty} \lambda\cdot b\in F\right\}\,.\end{equation}
\end{remark}

\section{Fixed points for the subgroups of $S^1$}
\label{podgrupy}
  Suppose $X$ is a smooth complete complex algebraic variety. In \cite{B-B2} there were studied the fixed point sets of the subgroups of $\T$. Originally $X$ was assumed to be projective, but the results hold for complete varieties by \cite[Lemma 8]{Su}. Our goal is to show, that the information about the plus-decompositions of fixed point sets is encoded in the localized Hirzebruch classes.

 Let $ G\subset \T$ be the subgroup of order $k$. Then the fixed point set $X^G$ is again a $\T$-manifold, with $(X^G)^\T=X^\T$. The plus-decomposition of $X^G$ is given by \begin{equation}(X^G)^+_F=X^+_F\cap X^G\,.\end{equation} The dimension of the normal bundle of $F$ in $(X^G)_F^+$ is equal to
\begin{equation}\sum_{k|w,\; w>0} n_w(F)\,,\end{equation}
where $n_w(F)=\dim(\nu_F^w)$. In particular for $k=1$
\begin{equation}\sum_{w>0} n_w(F)=n^+(F)\,,\qquad \sum_{w<0} n_w(F)=n^-(F)\,.\end{equation}
Knowing all the numbers $n_w(F)$ and the cohomology of the fixed point set components one can compute (by Theorem \ref{BBdec}) the cohomology of $X^G$ for all $G\subset\T$. We claim that the numbers $n_w(F)$ are encoded in the localized Hirzebruch class. First let us note that applying the limit procedure
to
\begin{equation}\label{produkt} \frac{td_y^\T(X)_{|F}}{eu^\T(\nu_F)}=td_y(F)\cdot\frac{td_y(\nu(F))}{eu^\T(\nu_F)}\end{equation}
    we obtain $td_y(F)(-y)^{n^+(F)}$ (by Theorem \ref{gladloc}). Plugging $y=-1$ we obtain $td_y(F)$, the first factor of (\ref{produkt}). Therefore the expression (\ref{produkt}) contains information about its second factor
\begin{equation}\label{produktsec}
\frac{td_y(\nu(F))}{eu^\T(\nu_F)}\,.\end{equation}
The degree zero part of (\ref{produktsec}) with respect to the gradation in $H^*(F)$  is equal to
\begin{equation}\prod_{w\in\Z}\left(\frac{1+y\,e^{-w t}}{1-e^{-w t}}\right)^{n_w(F)}
=\prod_{w\in\Z}\left(\frac{\theta^w+y}{\theta^w-1}\right)^{n_w(F)}\,,\end{equation}
where $\theta=e^{t}$. Let us denote this expression by $\Delta({\underline n})$.
It is an exercise in calculus to show that the numbers $n_w$ are determined by $\Delta({\underline n})$. We just note the following equalities:
\begin{itemize}
\item The function $\Delta({\underline n})_{|y=0}$  has poles at roots of unity.
The order of the pole at a primitive root of order $k$ is equal to
$\sum_{w>0,\;k|w}(n_{-w}+n_w)\,,$
\item the coefficient at $\theta^k$, $k>0$ in $\frac{\partial}{\partial y}\Delta({\underline n})_{|y=-1}$ is equal to
$\sum_{w>0,\;w|k}(n_{-w}-n_w)\,.$
\end{itemize}
From these data the numbers $n_{-w}- n_{w}$ and $n_{-w}+ n_{w}$ can be computed inductively with respect to divisibility of $w$. Finally all the exponents $n_w$ may be found.

\begin{remark}\rm \label{kontr} The mapping
\begin{equation}\prod_{w\in\Z}\left(\frac{\theta^w+y}{\theta^w-1}\right)^{n_w(F)}\quad\mapsto\quad \prod_{w>0,\,k|w} (-y)^{n_w}\end{equation}
cannot be extended to a linear functional. There are relations, such as
\begin{equation}\Delta(1,4)-\Delta(1,3)+\Delta(2,2)+\Delta(3,4)-2\,\Delta(2,4)=0.
\end{equation}
The above combination
for $k=2$ would have to be sent to
\begin{equation}-y-1+y^2-y-2y^2\neq 0.\end{equation}
Therefore
 for singular varieties it is not possible to compute the classes of the cells $(X^G)^+_F=(X^+_F)^{G}$ from the localized Hirzebruch class. More specifically, there is no linear function sending
\begin{equation}\frac{td^\T_y(X)_{|F}}{eu^\T(\nu_F)}=td_y(F)\prod_i\frac{1+y\,\theta^{-w_i}e^{-x_i}}{1-\theta^{-w_i}e^{-x_i}}\,.\end{equation}
 to
 \begin{equation}td_y((X^G)^+_F\to F)=td_y(F) \prod_{w>0,\,k|w} (-y)^{n_w}.\end{equation}
\end{remark}

\section{Beyond algebraic geometry}

We would like to mention a decomposition of $\chi_y$-genus in the case of symplectic manifolds. We do not consider singular spaces, since for such spaces it is not clear how to extend the definition of $\chi_y$-genus. We will just discuss the formula (\ref{hirdec}).

Let $X$ be a symplectic manifold with a Hamiltonian action of $T=S^1$ with a moment map (i.e. the Hamiltonian) $\mu$. We fix a Riemannian metric on $X$ which agrees with the symplectic structure. This way the tangent bundle to $X$ is equipped with a complex structure. We can assume, averaging if necessary, that the scalar product is invariant with respect to the action of $T$.
Now it makes perfect sense to apply the Hirzebruch class $td_y(-)$ to the tangent bundle $TX$ and compute $\chi_y(X)$. The equivariant version $td_y^T(-)$ is available as well. Theorem \ref{locth} allows to compute $\chi_y(X)$ as a sum of local contributions.

Let $F$ be a component of $X^T$. For $x\in X^T$ the tangent space $T_x X$ has a complex structure and as in the case of algebraic manifolds we define the number $n^+(F)$ as the dimension of the subrepresentation with positive weights. The number $n^+(F)$ can be computed in another way, using the gradient flow.
Consider the flow $f_t:X\to X$, $t\in \R$ associated to the vector field $grad(\mu)$. Let
\begin{equation}X^+_F=\{x\in X\,|\,\lim_{t\to-\infty} f_t(x)\in F\}\end{equation}
be the unstable subset. Then
\begin{equation}2n^+(F)=\dim (X^+_F)-\dim(F)\,.\end{equation}
It depends only on the Hessian of $\mu$ restricted to a normal slice of $F$. This restriction of the Hessian is nondegenerate due to \cite[Theorem 32.6]{GuSt}.
Applying Theorem \ref{locth} and the formula for the limit of the localized Hirzebruch class (Corollary \ref{granice}) we obtain the decomposition
\begin{equation}\chi_y(X)=\sum_{F\text{ \rm component of }X^T} \chi_y(F)\cdot (-y)^{n^+(F)}\,.\end{equation}
\begin{remark}\rm

The corresponding formula for Betti numbers
\begin{equation}b_k(X)=\sum_{F\text{ \rm component of }X^T} b_{k-2n^+(F)}(F)\end{equation} was proved by
Frankel \cite{Fr} for K\"ahler manifolds (with $b_1(X)=0$ or under the assumption that $X^T\not=\emptyset$), but his proof works for symplectic manifolds with Hamiltonian action. Frankel treats the moment map $\mu$ as a Morse function.
The critical set of $\mu$ is equal to the fixed point set $X^T$ and consists of nondegenerate critical submanifolds in the sense of Bott \cite{Bt}.
The proof of Frankel relies on Morse-Bott inequalities
\begin{equation}b_k(X)\leq \sum_{F\text{ \rm component of }X^T} b_{k-2n^+(F)}(F)\end{equation}
and opposite inequalities which hold in general topological context.
 The exact analogue of the Bia\l ynicki-Birula decomposition in the case of Hamiltonian actions is mentioned in \cite[Theorem 32.5]{GuSt} recalling Bott's work. A relevant reference for this circle of problems is the book \cite{Ki}, especially \S I.5.
\end{remark}

Further weakening of assumptions on $X$, i.e. assuming only that $X$ is a smooth oriented manifold, leads to the formula for signature ($y=1$) and Euler characteristic ($y=-1$).  Taking the limit when $t\to -\infty$ we obtain:
\begin{equation}\sigma(X)=\sum_{F\text{ \rm component of }X^T} \pm\sigma(F)\,,\end{equation}
with sign\footnote{Orientations of the ambient manifold $M$ and the submanifold $F$ induce an orientation of the normal bundle. The action of $S^1$ on the normal bundle induces another orientation. The sign depends on whether these orientations agree.}  depending on the character of the action of $S^1$ and the orientation of the normal bundle to $F$. This equality is the limit of \cite[Formula (45)]{Wi}. We note that localization of signature was studied in  \cite{AtSi} for finite order automorphisms, see Theorem 6.12 and its corollaries, in particular the formula 7.8 for circle action on a manifold with boundary. The equality $\sigma(X)=\sigma(X^T)$ is mentioned by Hirzebruch in the report for Mathematical Reviews.

Forgetting the orientation we are left only with the formula for Euler characteristic
\begin{equation}\chi(X)=\sum_{F\text{ \rm component of }X^T} \chi(F)\,,\end{equation}
which was probably the starting point for all localization theorems for torus actions.

%\bibliographystyle{alpha}
%\bibliography{hirbbbib}\end{document}

\end{document}